\documentclass{article}

\newtheorem{theorem}{Theorem}[section]
\newtheorem{corollary}[theorem]{Corollary}

\newtheorem{lemma}[theorem]{Lemma}
\newtheorem{example}[theorem]{Example}
\newtheorem{proposition}[theorem]{Proposition}
\newtheorem{remark}[theorem]{Remark}
\newtheorem{algorithm}[theorem]{Algorithm}

\def\qed{\hbox{}\nobreak [\kern-.4mm] \par \goodbreak \smallskip}

\title{{\bf Isomorphism classes of $A$-hypergeometric 
 systems}}
\author{Mutsumi Saito}

\date{December 27, 1999}

\begin{document}

\maketitle

\begin{abstract}
For a finite set $A$ of integral vectors,
Gel'fand, Kapranov and Zelevinskii defined a system of differential
equations with a parameter vector as a $D$-module,
which system is called an 
{\it $A$-hypergeometric} (or a {\it GKZ hypergeometric}) {\it system}.
Classifying the parameters according to the $D$-isomorphism classes of
their corresponding $A$-hypergeometric systems
is one of the most fundamental problems in the theory.
In this paper we give a combinatorial answer for the problem,
and illustrate it 
in two particularly simple cases:
the normal case and the monomial curve case.
\end{abstract}

\section{Introduction}

For a finite set $A$ of integral vectors,
Gel'fand, Kapranov and Zelevinskii defined a system of differential
equations with a parameter vector as a $D$-module,
which system is called an 
{\it $A$-hypergeometric} (or a {\it GKZ hypergeometric}) {\it system}
(\cite{GZK-def}).
Many authors studied $D$-invariants of the $A$-hypergeometric systems:
In Cohen-Macaulay case,
Gel'fand, Kapranov and Zelevinskii determined the characteristic cycles
(\cite{GZK})
and proved the irreducibility of
the monodromy representations for nonresonant parameters (\cite{GKZ});
Adolphson proved the rank of an $A$-hypergeometric system
equals the volume of the convex hull of $A$ in the
semi-nonresonant case (\cite{Adolphson});
The author, Sturmfels and Takayama 
scrutinized the ranks in \cite{sst-book};
Cattani, D'Andrea, and Dickenstein determined rational solutions and
algebraic solutions in monomial curve case (\cite{CDD}), and
recently Cattani, Dickenstein, and 
Sturmfels in \cite{CDS} considered when an $A$-hypergeometric system 
has a rational solution other than Laurent polynomial solutions.

The purpose of this paper is
to classify $A$-hypergeometric systems with respect to $D$-isomorphisms.
This is one of the most fundamental problems in the theory.
Under the assumption that the finite set $A$ lies in a hyperplane off the
origin,
we shall give a combinatorial answer for this problem, and illustrate it 
in two particularly simple cases:
the normal case and the monomial curve case.

Throughout the paper, we consider the finite set $A$ fixed.
In Section 2, we define a finite set $E_\tau(\beta)$ for
a parameter $\beta$ and a face $\tau$ of the cone generated by $A$.
Then our main theorem (Theorem \ref{theorem:main})
states that two $A$-hypergeometric systems
corresponding to parameters $\beta$ and $\beta'$
are $D$-isomorphic if and only if
$E_\tau(\beta)$ equals $E_\tau(\beta')$ for all faces $\tau$.
In Section 2, we prove the only-if-part of the theorem
and state some basic properties of the set $E_\tau(\beta)$.

Sections 3 and 4 are devoted to the study of
the algebra of contiguity operators, which algebra is called the {\it
symmetry algebra}.
In Section 3, we summarize some known facts on the symmetry algebra.
We introduce the {\it $b$-ideals} in Section 4 and prove
their elements correspond to contiguity operators.
Furthermore we describe each $b$-ideal in terms of the
standard pairs of a certain monomial ideal.
Using this description,
we give the proof of the if-part of our main theorem
in the end of Section 4.

In Sections 5 and 6, we illustrate our main theorem in the normal case
and the monomial curve case respectively,
since the theorem reduces to relatively simple forms in both cases.

The author is very grateful to Professor Nobuki Takayama for 
his fruitful comments.

\section{Main theorem}

We work over a field ${\bf k}$ of characteristic zero.
Let $A=(a_1,\ldots, a_n)=(a_{ij})$ be an integer $d\times n$-matrix
of rank $d$. 
We assume that all $a_j$ belong to one hyperplane off the origin 
in ${\bf Q}^d$.
We denote by $I_A$ the toric ideal in
${\bf k}[\partial]={\bf k}[\partial_1,\ldots,\partial_n]$, that is
$$
I_A=\langle \partial^u-\partial^v\, |\, Au=Av, \, u, v\in {\bf N}^n
\rangle \subset {\bf k}[\partial].
$$
For a column vector $\beta={}^t (\beta_1,\ldots,\beta_d)\in {\bf k}^d$,
let $H_A(\beta)$ denote the left ideal of the Weyl algebra
$$
D={\bf k}\langle x_1,\ldots, x_n,\partial_1,\ldots,
\partial_n\rangle
$$
generated by $I_A$ and $\sum_{j=1}^n a_{ij}\theta_{j} -\beta_i$
($i=1,\ldots, d$)
where $\theta_j =x_j\partial_j$.
The quotient $M_A(\beta)=D/H_A(\beta)$ is called
the {\it $A$-hypergeometric system with parameter $\beta$}.

We denote the set $\{\, a_1,\ldots, a_n\,\}$ by $A$ as well.
Let $\tau$ be a face of the cone 
\begin{equation}
{\bf Q}_{\geq 0}A=
\{\, \sum_{j=1}^n c_ja_j\, |\, c_j\in {\bf Q}_{\geq 0}\,\}.
\end{equation}
For a parameter $\beta\in {\bf k}^d$,
we consider the following set:
\begin{equation}
E_\tau (\beta):=
\{\, \lambda \in {\bf k}(A\cap\tau) /{\bf Z}(A\cap\tau)\, |\,
\beta-\lambda \in {\bf N}A+{\bf Z}(A\cap\tau)\,\}.
\end{equation}
Here 
${\bf N}=\{\, 0, 1, 2,\ldots\,\}$ and
we agree that ${\bf k}(A\cap\tau)={\bf Z}(A\cap\tau)=\{\, 0\,\}$
when $\tau=\{\, 0\,\}$.

The following is the main theorem in this paper.

\begin{theorem}
\label{theorem:main}
The $A$-hypergeometric systems $M_A(\beta)$ and $M_A(\beta')$ are
isomorphic as $D$-modules if and only if
$E_\tau(\beta)=E_\tau(\beta')$ for all faces $\tau$ of the cone
${\bf Q}_{\geq 0}A$.
\end{theorem}

Before the proof, we recall the formal series solutions $\phi_v$
defined in \cite{sst-book}.
For $v\in {\bf k}^n$, its {\it negative support} 
${\rm nsupp}(v)$ is the set of indices
$i$ with $v_i\in {\bf Z}_{<0}$.
When ${\rm nsupp}(v)$ is minimal with respect to inclusions
among ${\rm nsupp}(v+u)$ with $u\in {\bf Z}^n$ and $Au=0$,
$v$ is said to have {\it minimal negative support}.
For $v$ with minimal negative support,
we define a formal series
\begin{equation}
\phi_v=
\sum_{u\in N_v}\frac{[v]_{u_-}}{[v+u]_{u_+}}x^{v+u}.
\end{equation}
Here 
$$
N_v=\{\, u\in {\bf Z}^n\, |\,
Au=0,\, {\rm nsupp}(v)={\rm nsupp}(v+u)\,\},
$$
and
$u_+, u_-\in {\bf N}^n$ satisfy $u=u_+ -u_-$ with disjoint supports,
and
$[v]_w=\prod_{j=1}^n v_j(v_j-1)\cdots (v_j-w_j+1)$
for $w\in {\bf N}^n$.
Proposition 3.4.13 of \cite{sst-book}
states that
the series $\phi_v$ is a formal solution of $M_A(Av)$.

{\bf Proof.}
Here we prove the only-if-part of the theorem.
The proof of the if-part will be given in the end of
Section \ref{section:b-ideals}.

We suppose
that $\lambda\in E_\tau(\beta)\setminus E_\tau(\beta')$,
and we shall prove $M_A(\beta)$ and $M_A(\beta')$ are not isomorphic.

Represent $\lambda$ as $\sum_{a_j\in \tau}l_j a_j$.
Consider the direct product
$$
R_{\tau, \lambda}:=
\prod_{u\in {\bf Z}^n,\, u_j \in {\bf N} \,
(a_j\notin \tau)}
{\bf k}x^{l +u}.
$$
Here we put $l_j=0$ for $a_j\notin \tau$.
Note that $R_{\tau, \lambda}$ has the natural $D$-module
structure.
There exists $u\in {\bf Z}^n$ with $u_j\in {\bf N}$ 
($a_j\notin \tau$) such that
$\beta =A(l+u)$ and $l+u$ has minimal negative support.
Then the series $\phi_{l+u}\in R_{\tau, \lambda}$
is a formal solution of $M_A(\beta)$, and hence
${\rm Hom}_D(M_A(\beta), R_{\tau, \lambda})\not= 0$.
On the other hand,
${\rm Hom}_D(M_A(\beta'), R_{\tau, \lambda})= 0$
since $A(l+u)\not= \beta'$
for any $u\in {\bf Z}^n$ with $u_j \in {\bf N}$ $(a_j\notin \tau)$.
Therefore $M_A(\beta)$ and $M_A(\beta')$ are not isomorphic.
\qed

In the remainder of this section, we collect some properties
of the set $E_\tau(\beta)$.
We call a face of ${\bf Q}_{\geq 0}A$ of dimension $d-1$, a facet.
Recall that for a facet $\sigma$
the linear form $F_\sigma$ satisfying the following conditions is 
unique and called the {\it primitive integral support function}:

\begin{enumerate}
\item $F_\sigma ({\bf Z}A)={\bf Z}$, 
\item $F_\sigma(a_j)\geq 0$ for all $j=1,\ldots, n$,
\item $F_\sigma (a_j)=0$ for all $a_j\in \sigma$.
\end{enumerate}

\begin{proposition}
\label{remark:whole-cone}
\begin{enumerate}
\item
Each
$E_{{\bf Q}_{\geq 0}A}(\beta)$ consists of one element.
The equality $E_{{\bf Q}_{\geq 0}A}(\beta)=
E_{{\bf Q}_{\geq 0}A}(\beta')$ means $\beta -\beta'\in {\bf Z}A$.
\item
$E_{\{\, 0\,\}}(\beta)=\{\, 0\,\}$ or $\emptyset$.
$E_{\{\, 0\,\}}(\beta)=\{\, 0\,\}$ if and only if $\beta\in {\bf N}A$.
\item
For a facet $\sigma$, $E_\sigma(\beta)\not=\emptyset$
if and only if $F_\sigma (\beta)\in F_\sigma({\bf N}A)$.
\item
For faces $\tau\subset\sigma$,
there exists a natural map from $E_\tau(\beta)$ to $E_\sigma(\beta)$.
In particular, if $E_\tau(\beta)\not=\emptyset$, then
$E_\sigma(\beta)\not=\emptyset$.
\item
For any $\chi\in {\bf N}A$, there exists a natural 
inclusion from $E_\tau(\beta)$ to $E_\tau(\beta+\chi)$.
\end{enumerate}
\end{proposition}

{\bf Proof.}
All statements follow directly from the definition of $E_\tau(\beta)$.
\qed

\begin{proposition}
\label{remark:index-1}
\begin{enumerate}
\item
\begin{equation}
|E_\tau(\beta)|\leq
[({\bf Q}(A\cap\tau))\cap {\bf Z}A : {\bf Z}(A\cap \tau)].
\end{equation}
\item
Assume $({\bf Q}(A\cap\tau))\cap {\bf Z}A  ={\bf Z}(A\cap \tau)$.
If $\beta-\beta'\in {\bf Z}A$, and if neither $E_\tau(\beta)$ 
nor $E_\tau(\beta')$ is empty,
then $E_\tau(\beta)=E_\tau(\beta')$.
\end{enumerate}
\end{proposition}

{\bf Proof.}
\begin{enumerate}
\item
Let $\lambda,\lambda'\in E_\tau(\beta)$.
Then $\lambda-\lambda'\in ({\bf k}(A\cap\tau))\cap {\bf Z}A$.
By Cram\'er's formula, $({\bf k}(A\cap\tau))\cap {\bf Z}A
=
({\bf Q}(A\cap\tau))\cap {\bf Z}A$.
\item
Let $E_\tau(\beta)=\{\, \lambda\,\}$, $E_\tau(\beta')=\{\, \lambda'\,\}$.
Since $\beta-\beta'\in {\bf Z}A$, there exist
$\chi,\chi'\in {\bf N}A$ such that
$\beta+\chi =\beta' +\chi'$.
Then $\{\, \lambda\,\} = E_\tau(\beta+\chi)= E_\tau(\beta'+\chi')
=\{\, \lambda'\,\}$
by Proposition \ref{remark:whole-cone} (5).
\end{enumerate}
\qed

\begin{example}
Let 
$$
A=\left(
\begin{array}{cccc}
1 & 1 & 1 & 1\\
0 & 0 & 1 & 2\\
0 & 1 & 1 & 0
\end{array}
\right).
$$
There are four facets:
\begin{eqnarray}
\sigma_{12}:&=& {\bf Q}_{\geq 0}a_1+{\bf Q}_{\geq 0}a_2,\\
\sigma_{23}:&=& {\bf Q}_{\geq 0}a_2+{\bf Q}_{\geq 0}a_3,\\
\sigma_{34}:&=& {\bf Q}_{\geq 0}a_3+{\bf Q}_{\geq 0}a_4,\\
\sigma_{14}:&=& {\bf Q}_{\geq 0}a_1+{\bf Q}_{\geq 0}a_4,
\end{eqnarray}
and four one-dimensional faces: ${\bf Q}_{\geq 0}a_1,\ldots, 
{\bf Q}_{\geq 0}a_4$.
For all faces $\tau$ but $\sigma_{14}$,
the indices $[({\bf Q}(A\cap\tau))\cap {\bf Z}A : {\bf Z}(A\cap\tau)]$
are one.
Hence for $\beta\in {\bf N}A$, $E_\tau(\beta)=\{\, 0\,\}$
for all faces $\tau\not= \sigma_{14}$.
The quotient
$({\bf Q}(A\cap\sigma_{14}))\cap {\bf Z}A / {\bf Z}(A\cap\sigma_{14})$
has two elements and can be represented by $0$ and ${}^t (1,1,0)$.
Since $a_2-{}^t (1,1,0)=a_3 -a_4$, and $a_3-{}^t (1,1,0)=a_2-a_1$,
we obtain 
$E_{\sigma_{14}}(a_2)=E_{\sigma_{14}}(a_3)=\{\, 0, {}^t (1,1,0)\,\}$.
Proposition \ref{remark:whole-cone} (5) implies that for
$\beta\in {\bf N}A$,
$E_{\sigma_{14}}(\beta)=\{\, 0\,\}$ if and only if
$\beta\in {\bf N}a_1 +{\bf N}a_4$,
 otherwise $E_{\sigma_{14}}(\beta)=\{\, 0, {}^t (1,1,0)\,\}$.
Therefore ${\bf N}A$ splits into two isomorphism classes in this case.
\end{example}

Recall that a parameter $\beta$ is said to be {\it
nonresonant} (respectively {\it semi-nonresonant})
if $\beta\notin {\bf Z}A +{\bf k}(A\cap\sigma)$
(respectively
$\beta\notin ({\bf Z}A\cap {\bf Q}_{\geq 0}A) +{\bf k}(A\cap\sigma)$
) for any facet $\sigma$,
or equivalently,
if
$F_\sigma(\beta)\notin {\bf Z}$ (respectively
$F_\sigma(\beta)\notin {\bf N}$
) for any facet $\sigma$.
Hence the nonresonance implies the semi-nonresonance.

\begin{proposition}
\label{proposition:semi-nonresonance}
If $\beta$ is semi-nonresonant, then
$E_\tau(\beta)=\emptyset$ for all proper faces $\tau$
of ${\bf Q}_{\geq 0}A$.
\end{proposition}

{\bf Proof.}
The semi-nonresonace clearly implies
$E_\sigma(\beta)=\emptyset$ for all facets $\sigma$.
Proposition \ref{remark:whole-cone} (4) finishes the proof.
\qed

\begin{corollary}
\label{corollary:semi-nonresonance}
Let $\beta$ and $\beta'$ be semi-nonresonant.
Then $M_A(\beta)$ and $M_A(\beta')$ are isomorphic
if and only if $\beta-\beta'\in {\bf Z}A$.
\end{corollary}

\begin{proposition}
If a parameter $\beta$ satisfies
\begin{equation}
E_\tau(\beta)=\emptyset \quad \mbox{for all proper faces $\tau$,}
\end{equation}
then 
\begin{enumerate}
\item
for any $\chi\in {\bf N}A$,
$M_A(\beta-\chi)$ is isomorphic to $M_A(\beta)$.
\item
Recall that all elements of $A$ lie on one hyperplane
$H$ off the origin.
We normalize the volume of a polytope on $H$
so that a simplex whose vertices affinely span the lattice
$H\cap {\bf Z}A$ has volume one.
Then the rank of $M_A(\beta)$, i.e. the rank of the solution sheaf of
$M_A(\beta)$, equals the volume of the convex hull of $A$.
\end{enumerate}
\end{proposition}

{\bf Proof.}
(1)\quad
By Proposition \ref{remark:whole-cone} (5),
$E_\tau(\beta-\chi)=\emptyset$ for all proper faces $\tau$.
Hence by Proposition \ref{remark:whole-cone} (1),
we deduce the statement from Theorem \ref{theorem:main}.

The proof of (2) is the same as that of Theorem 4.5.2
of \cite{sst-book} (p. 185).
\qed

\section{Symmetry algebra}

We consider
the algebra of contiguity operators,
which algebra is called the symmetry algebra.
It controls isomorphisms among
$A$-hypergeometric
systems with different parameters.
We have investigated
the symmetry algebra of normal $A$-hypergeometric
systems in \cite{SymAlg}.
The proofs of some results in \cite{SymAlg}
remain valid without the normality condition.
In this section, we summarize such results.

Let
$$
\tilde{S}:=\{\, P\in D\, |
\, I_A P \subset D I_A\, \}.
$$
Then $\tilde{S}$ is an associative algebra
and $\tilde{S}\cap D I_A$ is its two-sided ideal.
We call $S:=\tilde{S}/\tilde{S}\cap D I_A$
the {\it symmetry algebra} of $A$-hypergeometric
systems.
The symmetry algebra $S$ is nothing but the algebra
${\rm End}_D(D/DI_A)$.
We remark that $D/DI_A$ can be considered as the system of
differential equations for the generating functions of
$A$-hypergeometric functions.

In what follows, we denote simply by $P$,
the element of $D/DI_A$ represented by
$P\in D$.
For $\chi\in {\bf N}A$,
all $\partial^u$ with $Au=\chi$ represent the same element
in $D/DI_A$. Hence we sometimes denote it by $\partial^\chi$.

\begin{proposition}
\label{corollary:1.3}
\begin{enumerate}
\item
$\partial_1,\ldots,\partial_n\in S$.
\item
$\sum_{j=1}^n a_{ij}\theta_j\in S$ for all $i=1,\ldots, d$.
\item
The morphism from the polynomial ring 
${\bf k}[\, s\,]={\bf k}[\, s_1,\ldots, s_d\,]$ to $S$ mapping $s_i$ to
$\sum_{j=1}^n a_{ij}\theta_j$  ($i=1,\ldots,d$)
is injective.
\end{enumerate}
\end{proposition}

{\bf Proof.}
See Lemma 1.1 in \cite{SymAlg} for (1) and (2), and
Corollary 1.3 in \cite{SymAlg} for (3).
\qed

We consider that ${\bf Z}A$ is the character group
of the algebraic torus $T=\{\, (t_1,\ldots, t_d)\, |\, 
t_1,\ldots, t_d\in {\bf k}^\times\,\}$.
Let $N$ be the dual group of ${\bf Z}A$, and $s_1,\ldots, s_d$
the basis of ${\bf k}\otimes_{\bf Z} N$
dual to the standard basis of ${\bf k}^d={\bf k}\otimes_{\bf Z} {\bf Z}A$.
Under the identification of ${\bf k}\otimes_{\bf Z} N$ with
the Lie algebra of $T$ (\cite{Oda}),
each $s_i$ equals $t_i\frac{\partial}{\partial t_i}$.
The morphism
in Proposition \ref{corollary:1.3} (3)
is induced from the differential of the injective morphism:
\begin{equation}
T\ni t\longmapsto (t^{a_1},\ldots, t^{a_n})\in ({\bf k}^\times)^n.
\end{equation}
We thus consider ${\bf k}[\, s\,]$ as a subspace of $S$ and, accordingly, 
as a subspace of $D/DI_A$.
For each $\chi\in {\bf Z}A$, 
we define the weight space $S_{\chi}$ with 
weight $\chi$ by
$$
S_\chi :=\{\, P\in S\, |\, [s,P]=\chi(s)P\quad (\forall s\in N)\,\}.
$$
Here the bracket $[P,Q]$ means $PQ-QP$.

\begin{remark}
Note that the multiplication by $P\in S_\chi$ from the right
defines a $D$-homomorphism from $M_A(\beta+\chi)$ to
$M_A(\beta)$. Hence $P(\psi_\beta)$ is a solution of $M_A(\beta+\chi)$
for a solution $\psi_\beta$ of $M_A(\beta)$.
In this sense, the operator $P$ is a contiguity operator shifting
parameters by $\chi$.
\end{remark}

\begin{theorem}
\label{proposition:prop2-9}
\begin{enumerate}
\item
The symmetry algebra $S$ has no zero-divisors.

\item
The symmetry algebra $S$ has the following weight space decomposition:
\begin{equation}
S = \bigoplus_{\chi\in {\bf Z}A}S_\chi.
\end{equation}

\item
The weight space
$S_0$ equals the polynomial ring ${\bf k}[\, s\,]$.

\item
For each $\chi\in {\bf N}A$,
the weight space
$S_{-\chi}$ equals ${\bf k}[\, s\, ]\partial^\chi$.

\end{enumerate}
\end{theorem}

{\bf Proof.}
See Lemma 1.4, and Propositions 2.3, 2.4, 2.9 in \cite{SymAlg}.
\qed

The following proposition will be used in the next section.

\begin{proposition}[{\rm Proposition 2.6 in \cite{SymAlg}}]
\label{proposition:2.6}
The natural morphism
$$
D/DI_A \longrightarrow
{\bf k}\langle\, x, \partial^{\pm}\,\rangle/
{\bf k}\langle\, x,\partial^{\pm}\,\rangle I_A
$$
is injective where
${\bf k}\langle\, x, \partial^{\pm}\,\rangle$
is the algebra generated by $D$ and 
 $\partial_1^{-1},\ldots,\partial_n^{-1}$
with relations
$[\, x_i, \partial_j^{-1}\,]=\delta_{ij}\partial_j^{-2}$
$(i,j=1,\ldots, n)$.
\end{proposition}

\section{$b$-Ideals}
\label{section:b-ideals}

We have seen in Theorem \ref{proposition:prop2-9} that
the symmetry algebra $S$ has a weight decomposition
with respect to ${\bf Z}A$, and that
each $S_\chi$ for $-\chi \in {\bf N}A$
is the free ${\bf k}[s]$-module of rank one 
with basis $\partial^{-\chi}$.
Next we wish to compute the weight space $S_\chi$ for arbitrary $\chi$.
Suppose that $E\in S_\chi$ and $\chi=\chi_+-\chi_-$ with $\chi_+, \chi_-
\in {\bf N}A$.
Then the operator $E\partial^{\chi_+}$ belongs to $S_{-\chi_-}$.
Hence 
by Theorem \ref{proposition:prop2-9} (4),
there exists a polynomial $b\in {\bf k}[s]$
such that $E\partial^{\chi_+}= b\partial^{\chi_-}$.
Such polynomials $b$ varying $E\in S_\chi$ form an ideal
of ${\bf k}[s]$.
We shall define the {\it $b$-ideal} $B_\chi$ below to be such an ideal.

Fix any $\chi\in {\bf Z}A$, and
define an ideal $I_\chi$ of ${\bf k}[\partial]$ by

\begin{equation}
I_\chi:=I_A+M_\chi
\end{equation}
where
\begin{equation}
M_\chi:=\langle \,
\partial^u\, |\,
Au\in \chi+{\bf N}A
\,\rangle.
\end{equation}
Define the ideal $B_\chi$ of {\it $b$-polynomials} by

\begin{equation}
B_\chi:= {\bf k}[s]\cap DI_\chi.
\end{equation}

\begin{proposition}
\label{proposition:EandB}
Let $\chi=\chi_+ -\chi_-$ with $\chi_+,\chi_-\in {\bf N}A$.
For $b\in B_\chi$, there exists a unique operator $E\in S_\chi$ 
such that $b\partial^{\chi_-} = E \partial^{\chi_+}$.
The operator $E$ is independent of the expression $\chi=\chi_+ -\chi_-$.

Moreover any operator in $S_\chi$ can be obtained in this way.
\end{proposition}

{\bf Proof.}
Since $b\partial^{\chi_-}\in DI_\chi\partial^{\chi_-}
\subset DI_A + D \partial^{\chi_+}$,
there exists an operator $E\in D$ such that 
$b\partial^{\chi_-} = E \partial^{\chi_+}$.
The uniqueness,
the independence,
and $E\in S_\chi$ follow from the equality $E=b\partial^{\chi}$ in
${\bf k}\langle x, \partial^\pm \rangle$
and Proposition \ref{proposition:2.6}.

Let $E\in S_\chi$ and $\chi=\chi_+-\chi_-$ with $\chi_+,\chi_-\in {\bf N}A$.
Then $E\partial^{\chi_+}\in S_{-\chi_-}$.
By Theorem \ref{proposition:prop2-9} (4),
there exists a polynomial $b\in {\bf k}[s]$
such that $E\partial^{\chi_+}=b\partial^{\chi_-}$.
Then $b\in I_\chi$ and thus $b\in B_\chi$.
\qed

We have the following algorithm of obtaining the operator
$E\in S_\chi$ corresponding to $b\in B_\chi$,
which generalizes Algorithm 3.4 in \cite{sst}.

\begin{algorithm}
Let $\chi=Au -Av$ and $u, v\in {\bf N}^n$.

Input: a polynomial $b\in B_\chi$.

Output: an operator $E\in S_\chi$ with $E\partial^{u}=
b\partial^{v}$.

\begin{enumerate}
\item
For $i=1,\ldots, n$,
compute a Gr\"obner basis ${\cal G}_i$ of $I_A$ with respect to
any reverse lexicographic term order with lowest variable $\partial_i$.
\item
Expand $b(\sum_{j}a_{1j}\theta_j,\ldots,\sum_{j}a_{dj}\theta_j)
\partial^v$
in ${\bf Q}\langle x,\partial\rangle$
into a ${\bf Q}$-linear combination of monomials
$x^l\partial^m$.
\item
$i:=1$, $E:=\mbox{the output of Step 2}$.

While $i\leq n$, do
\begin{enumerate}
\item
Reduce $E$ modulo ${\cal G}_i$ 
in ${\bf Q}\langle x,\partial\rangle$.
\item
The output of Step 3-(a) has $\partial_i^{u_i}$ as a right factor.
Divide it by $\partial_i^{u_i}$.
\item
$i:=i+1$, $E:=\mbox{the output of Step 3-(b)}$.
\end{enumerate}
\end{enumerate}
\end{algorithm}

The proof of the correctness is completely
analogous to that of Algorithm 3.4 in
\cite{sst}.

We thus reduce the study of $S_\chi$ to that of $B_\chi$,
and for the study of $B_\chi={\bf k}[s]\cap DI_\chi$, 
we study ${\bf k}[\theta]\cap DI_\chi$ first.
Since $M_\chi$ is the largest monomial ideal in $I_\chi$,
we have by Lemma 4.4.4 in \cite{sst-book},

\begin{proposition}
\label{proposition:lemma4-4-4}
\begin{equation}
{\bf k}[\theta]\cap DI_\chi =\widetilde{M_\chi}
\end{equation}
where $\widetilde{M_\chi}$ is the distraction of $M_\chi$,
i.e.,
$\widetilde{M_\chi}={\bf k}[\theta]\cap D M_\chi$.
\end{proposition}

For the study of $\widetilde{M_\chi}$, 
we recall the standard pairs of a monomial ideal.
Let $M$ be a monomial ideal of ${\bf k}[\partial]$.
Then a pair $(u, \tau)$ with $u\in {\bf N}^n$ and $\tau\subset \{\, 1,
\ldots, n\}$ is called a {\it standard pair} of $M$
if it satisfies the following conditions:
\begin{enumerate}
\item
$u_j=0$ for all $j\in \tau$. (We abbreviate this to $u\in {\bf N}^{\tau^c}$,
where ${}^c$ stands for taking the complement.)
\item
There exists no $v\in {\bf N}^\tau$ such that $\partial^{u+v}\in M$.
\item
For each $j\notin \tau$,
there exists $v\in {\bf N}^{\tau\cup\{ j\}}$ such that 
$\partial^{u+v}\in M$.
\end{enumerate}

For an algorithm of obtaining the set of standard pairs,
see \cite{Hosten-Thomas}.
Let ${\cal S}(M_\chi)$ denote the set of standard pairs of $M_\chi$.
By Corollary 3.2.3 in \cite{sst-book},
the distraction $\widetilde{M_\chi}$
is described as follows:

\begin{equation}
\label{equation:distraction}
\widetilde{M_\chi}=\bigcap_{(u,\tau)\in {\cal S}(M_\chi)}
\langle \,
\theta_i-u_i\, |\, i\notin \tau
\,\rangle.
\end{equation}

\begin{lemma}
\label{lemma:face}
Let $(u,\tau)$ be a standard pair of $M_\chi$.
Then $A{\bf Q}_{\geq 0}^\tau:=\sum_{j\in\tau}{\bf Q}_{\geq 0}a_j$
is a proper face of ${\bf Q}_{\geq 0}A$,
and moreover
$\tau =\{\, i\, |\, a_i\in A{\bf Q}_{\geq 0}^\tau\,\}$.
\end{lemma}

{\bf Proof.}
Suppose that $A{\bf Q}_{\geq 0}^\tau$ is not contained 
in any facet of ${\bf Q}_{\geq 0}A$.
Then there exists $\gamma\in A{\bf N}^\tau:=\sum_{j\in \tau}{\bf N}a_j$ 
such that
$F_\sigma(\gamma)>0$ for all facets $\sigma$.
Then $F_\sigma(Au +m\gamma)\gg 0$ for $m\gg 0$ and all facets $\sigma$.
By Lemma 1 in the appendix of \cite{Saito-Takayama},
$Au+m\gamma\in\chi+{\bf N}A$ for $m\gg 0$.
This contradicts $(u,\tau)$ is a standard pair of $M_\chi$.

Next we claim $(A{\bf Q}_{\geq 0}^{\tau^c})\cap (A{\bf Q}^\tau)
 =\{\, 0\,\}$,
which implies the lemma.
Suppose $(A{\bf Q}_{\geq 0}^{\tau^c})\cap (A{\bf Q}^\tau)
 \not=\{\, 0\,\}$.
Let $v\in {\bf N}^{\tau^c}$ be a nonzero element
satisfying $Av\in A{\bf Z}^\tau$.
Then there exists $w\in {\bf N}^\tau$ such that
$Aw\in Av+A{\bf N}^\tau$.
Since $A(u+mw)\notin M_\chi$ for any $m\in {\bf N}$,
$(Au +A{\bf N}^{\tau\cup \tau'})\cap M_\chi=\emptyset$
for $\tau'=\{ i\, |\, v_i\not= 0\}$.
This contradicts $(u,\tau)$ is a standard pair of $M_\chi$ again.
\qed

Thanks to Lemma \ref{lemma:face}, we regard the set $\tau$ 
of a standard pair $(u,\tau)$ as the 
proper face
$A{\bf Q}^\tau$ of ${\bf Q}_{\geq 0}A$.

For an ideal $I$ of ${\bf k}[s]$, we denote by $V(I)$ the zero set of $I$.
Proposition \ref{proposition:lemma4-4-4} and 
the equation (\ref{equation:distraction}) 
give the following 
prime decomposition of $B_\chi$ and irreducible decomposition of
the zero set $V(B_\chi)$.

\begin{theorem}
\label{theorem:b-ideal}
\begin{enumerate}
\item
\begin{equation}
B_\chi=
\bigcap_{(u,\tau)\in {\cal S}(M_\chi)}
\langle\, F_\sigma - F_\sigma(Au)\, |\, 
\mbox{$\sigma$ :facet $\supset\tau$}
\rangle.
\end{equation}
\item
\begin{equation}
V(B_\chi)=
\bigcup_{(u,\tau)\in {\cal S}(M_\chi)}
(Au +{\bf k}(A\cap\tau)).
\end{equation}
\end{enumerate}
\end{theorem}

{\bf Proof.}
From (\ref{equation:distraction}),
we only need to show
\begin{equation}
\label{equation:elimination}
{\bf k}[s]\cap
\langle \theta_i -u_i\, |\, i\notin\tau \rangle
=\langle F_\sigma -F_\sigma (Au)\, |\, \sigma\supset\tau
\rangle.
\end{equation}
First we have
\begin{equation}
V({\bf k}[s]\cap
\langle \theta_i -u_i\, |\, i\notin\tau \rangle)
=
Au+{\bf k}(A\cap \tau)
=V(\langle F_\sigma -F_\sigma (Au)\, |\, \sigma\supset\tau
\rangle).
\end{equation}
Hence
\begin{eqnarray}
{\bf k}[s]\cap
\langle \theta_i -u_i\, |\, i\notin\tau \rangle
&\supset&
\langle F_\sigma -F_\sigma (Au)\, |\, \sigma\supset\tau
\rangle
\nonumber\\
&=&
I(V(\langle F_\sigma -F_\sigma (Au)\, |\, \sigma\supset\tau
\rangle))
\nonumber\\
&=&
I(V({\bf k}[s]\cap
\langle \theta_i -u_i\, |\, i\notin\tau \rangle))
\end{eqnarray}
where $I$ stands for taking the defining ideal.
On the other hand, $J\subset I(V(J))$ is automatic for any ideal $J$.
We therefore obtain (\ref{equation:elimination}).
\qed

\begin{proposition}
\begin{enumerate}
\item
\begin{equation}
V(B_{\chi+\chi'})\subset V(B_\chi)\cup (V(B_{\chi'})+\chi)
\qquad \mbox{for $\chi,\chi'\in {\bf Z}A$}.
\end{equation}
\item
\begin{equation}
V(B_{\chi+\chi'})= V(B_\chi)\cup (V(B_{\chi'})+\chi)
\qquad \mbox{for $\chi,\chi'\in {\bf N}A$}.
\end{equation}
\end{enumerate}
\end{proposition}

{\bf Proof.}
\begin{enumerate}
\item
Let $p_\chi\in B_\chi$, $p_{\chi'}\in B_{\chi'}$,
and $P_\chi\in S_\chi$, $P_{\chi'}\in S_{\chi'}$ be
in the correspondence in Proposition \ref{proposition:EandB}.
Then
\begin{eqnarray}
P_{\chi}P_{\chi'} \partial^{\chi'_+}\partial^{\chi_+}
&=& P_{\chi} p_{\chi'}(s)\partial^{\chi'_-}\partial^{\chi_+}
\nonumber\\
&=& p_{\chi'}(s-\chi)P_{\chi}\partial^{\chi_+}\partial^{\chi'_-}
\nonumber\\
&=& p_{\chi'}(s-\chi) p_{\chi}(s)\partial^{\chi_-}\partial^{\chi'_-}.
\end{eqnarray}
Hence $p_{\chi'}(s-\chi) p_{\chi}(s)\in B_{\chi +\chi'}$.

\item
Let $p_{\chi+\chi'}\in B_{\chi+\chi'}$
and $P_{\chi+\chi'}\in S_{\chi+\chi'}$ be
in the correspondence in Proposition \ref{proposition:EandB}.
Then 
$$
p_{\chi+\chi'}=P_{\chi+\chi'}\partial^{\chi'}\cdot\partial^\chi.
$$
Hence $p_{\chi+\chi'}(s)\in B_\chi$.

Furthermore
$$
p_{\chi+\chi'}(s+\chi)\partial^\chi
=
\partial^\chi p_{\chi+\chi'}(s)=
\partial^\chi P_{\chi+\chi'}\partial^{\chi'}\partial^\chi.
$$
Hence $p_{\chi+\chi'}(s+\chi)=\partial^\chi P_{\chi+\chi'}\partial^{\chi'}$,
which impies $p_{\chi+\chi'}(s+\chi)\in B_{\chi'}$.
\end{enumerate}
\qed

\begin{proposition}
\label{proposition:4.10}
Let $\chi\in {\bf Z}A$.
Let $p_\chi\in B_\chi$, $p_{-\chi}\in B_{-\chi}$,
and $P_\chi\in S_\chi$, $P_{-\chi}\in S_{-\chi}$ be
in the correspondence in Proposition \ref{proposition:EandB}.
Then
\begin{equation}
P_{-\chi}P_\chi=p_\chi(s+\chi)p_{-\chi}(s).
\end{equation}
\end{proposition}

{\bf Proof.}
\begin{eqnarray}
P_{-\chi}P_\chi \partial^{\chi_+}
&=& P_{-\chi} p_{\chi}(s)\partial^{\chi_-}
\nonumber\\
&=& p_\chi(s+\chi)P_{-\chi}\partial^{\chi_-}
\nonumber\\
&=& p_\chi(s+\chi) p_{-\chi}(s)\partial^{\chi_+}.
\end{eqnarray}
Divide it by $\partial^{\chi_+}$ to obtain the conclusion.
\qed

For $\chi\in {\bf Z}A$,
define an ideal $B_{-\chi,\chi}$ by
\begin{equation}
B_{-\chi,\chi}:=
\langle \, p_\chi(s+\chi)p_{-\chi}(s)\, |\,
p_\chi\in B_\chi, p_{-\chi}\in B_{-\chi}\,\rangle.
\end{equation}
Then the following proposition is immediate from the definition
of $B_{-\chi,\chi}$.

\begin{proposition}
\label{corollary:4.11}
\begin{enumerate}
\item
\begin{equation}
V(B_{-\chi,\chi})=
(V(B_\chi)-\chi)\cup V(B_{-\chi}).
\end{equation}
\item
\begin{equation}
V(B_{-\chi,\chi})=
V(B_{\chi,-\chi})-\chi.
\end{equation}
\end{enumerate}
\end{proposition}

\begin{theorem}
\label{theorem:iso}
Let $\chi\in {\bf Z}A$.
If $\beta\notin V(B_{-\chi,\chi})$,
then two $A$-hypergeometric systems $M_A(\beta)$ and $M_A(\beta+\chi)$
are isomorphic.
\end{theorem}

{\bf Proof.}
First note that $\beta\notin V(B_{-\chi,\chi})$
is equivalent to $\beta+\chi\notin V(B_{\chi,-\chi})$
by Proposition \ref{corollary:4.11}.
Take polynomials $p_\chi\in B_\chi$ and $p_{-\chi}\in B_{-\chi}$
such that $p_{\chi}(\beta+\chi)p_{-\chi}(\beta)\not= 0$.
Let $P_\chi\in S_\chi$, $P_{-\chi}\in S_{-\chi}$ be
in the correspondence in Proposition \ref{proposition:EandB}.
Then by Proposition \ref{proposition:4.10},
we have the following equalities:
\begin{eqnarray}
P_{-\chi}P_\chi & = & p_\chi(s+\chi)p_{-\chi}(s),\\
P_{\chi}P_{-\chi} & = & p_{-\chi}(s-\chi)p_{\chi}(s).
\end{eqnarray}

The multiplications by $P_{-\chi}, P_\chi$ respectively
induce homomorphisms:

\begin{eqnarray}
f&:& M_A(\beta) \longrightarrow M_A(\beta+\chi),\\
g&:& M_A(\beta+\chi) \longrightarrow M_A(\beta).
\end{eqnarray}

Then 
\begin{equation}
g\circ f=p_{\chi}(\beta+\chi)p_{-\chi}(\beta)id_{M_A(\beta)}
\end{equation}
and
\begin{eqnarray}
f\circ g &=&
p_{-\chi}((\beta+\chi)-\chi)p_{\chi}(\beta+\chi)id_{M_A(\beta+\chi)}\\
&=&
p_{-\chi}(\beta)p_{\chi}(\beta+\chi)id_{M_A(\beta+\chi)}.
\end{eqnarray}
Hence $f$ and $g$ are isomorphisms.
\qed

Now we are ready to prove the if-part of our main theorem.

{\bf Proof of the if-part of Theorem \ref{theorem:main}.}

We suppose that $E_\tau(\beta)=E_\tau(\beta')$ for all faces.
Let $\chi:= \beta'-\beta$.
We claim $\beta\notin V(B_{-\chi})$.
Assume the contrary.
Then by Theorem \ref{theorem:b-ideal},
there exists a standard pair $(u,\tau)\in {\cal S}(M_{-\chi})$
such that $\beta- Au\in {\bf k}(A\cap \tau)$.
The equality $E_\tau(\beta)=E_\tau(\beta')$
implies that
there exists $v\in {\bf N}^n$ such that
$\beta -\beta'= A(u-v)$.
Hence the intersection
of $Au + {\bf N}(A\cap \tau)$ with $(\beta-\beta')+{\bf N}A$
is not empty.
This contradicts the standardness of $(u,\tau)$.
We have thus proved  $\beta\notin V(B_{-\chi})$.
By symmetry we have $\beta'\notin V(B_{\chi})$,
which is equivalent to $\beta\notin V(B_{\chi})-\chi$.
Hence $\beta\notin V(B_{-\chi,\chi})$
by Proposition \ref{corollary:4.11}.
From Theorem \ref{theorem:iso} we conclude
$M_A(\beta)$ is isomorphic to $M_A(\beta')$.
\qed

As a corollary of the proof of the if-part of Theorem \ref{theorem:main},
we obtain the following.

\begin{corollary}
If two $A$-hypergeometric systems $M_A(\beta)$ and $M_A(\beta')$ are
isomorphic, then there exists an operator $P\in S_{\beta'-\beta}$
such that the multiplication by $P$ from the right induces
an isomorphism from $M_A(\beta)$ to $M_A(\beta')$.
\end{corollary}

\section{Normal case}

In this section, we consider the normal case: 
\begin{equation}
\label{condition:normality}
{\bf N}A={\bf Z}A\cap {\bf Q}_{\geq 0}A.
\end{equation}

Many important examples are known to be normal,
such as Aomoto-Gel'fand systems, the $A$-hypergeometric
systems corresponding to ${}_{p+1}F_p$, Lauricella functions, etc.
(see \cite{saito-normal}, \cite{saito-fe}).
It will turn out below
that the parameter space can be classified in terms of
the primitive integral support functions $F_\sigma$
in the normal case.

\begin{lemma}
\label{lemma:normal}
Assume $A$ to be normal.
Then we have the following.
\begin{enumerate}
\item
$({\bf Q}(A\cap \tau))\cap {\bf Z}A$ equals ${\bf Z}(A\cap\tau)$
for all faces $\tau$.

\item
$F_\sigma({\bf N}A)={\bf N}$ for all facets $\sigma$.

\item
For a face $\tau$,
\begin{equation}
{\bf N}A +{\bf Z}(A\cap\tau)
=
{\bf Z}A\cap \bigcap_{\mbox{$\sigma$: facet $\supset\tau$}}
({\bf N}A +{\bf k}(A\cap\sigma)).
\end{equation}
\end{enumerate}
\end{lemma}

{\bf Proof.}
(1)\quad
Let $\chi\in ({\bf Q}(A\cap \tau))\cap {\bf Z}A$.
Add a vector $\chi'\in {\bf N}(A\cap\tau)$ to $\chi$
so that $\chi +\chi'\in {\bf Q}_{\geq 0}(A\cap \tau)$.
By the normality, we see $\chi+\chi'\in {\bf N}(A\cap\tau)$.
Hence $\chi$ belongs to ${\bf Z}(A\cap\tau)$.

(2)\quad
Let $\chi\in {\bf Z}A$ satisfy $F_\sigma(\chi)=1$.
For $\sigma'\not=\sigma$, there exists $a_j\in \sigma\setminus
\sigma'$. Hence there exists $\chi'\in {\bf N}(A\cap \sigma)$
such that $F_{\sigma'}(\chi+\chi')\geq 0$ for all facets $\sigma'$.
By the normality, $\chi +\chi'\in {\bf N}A$.
Since $F_\sigma(\chi+\chi')=1$, we obtain $F_\sigma({\bf N}A)={\bf N}$.

(3)\quad
Let $\chi\in {\bf Z}A$ satisfy $F_\sigma (\chi)\geq 0$ for
all facets containing $\tau$.
For a facet $\sigma$ not containing the face $\tau$,
 there exists $a_j\in \tau\setminus\sigma$.
Hence there exists a vector $\chi'\in {\bf N}(A\cap\tau)$
such that $F_\sigma(\chi +\chi')\geq 0$ for all facets $\sigma$
of the cone ${\bf Q}_{\geq 0}A$.
By the normality, $\chi +\chi'\in {\bf N}A$, and thus
$\chi\in {\bf N}(A\setminus A\cap\tau) +{\bf Z}(A\cap\tau)$.
\qed

\begin{theorem}
\label{theorem:normal}
Let $\beta, \beta'\in {\bf k}^d$. 
Then $M_A(\beta)\simeq M_A(\beta')$ if and only if
$\beta -\beta'\in {\bf Z}A$ and
$\{ \mbox{$\sigma$ : facet,} \, F_\sigma (\beta)\in {\bf N}\}=
\{ \mbox{$\sigma$ : facet,}\, F_\sigma (\beta')\in {\bf N}\}$.
\end{theorem}

{\bf Proof.}
By Proposition \ref{remark:whole-cone} (3), 
the only-if-part follows from Theorem \ref{theorem:main}.

Next we prove the if-part.
Suppose
$\beta -\beta'\in {\bf Z}A$ and
$\{ \mbox{$\sigma$ : facet,} \, F_\sigma (\beta)\in {\bf N}\}=
\{ \mbox{$\sigma$ : facet,}\, F_\sigma (\beta')\in {\bf N}\}$.
By Lemma \ref{lemma:normal} (1), (2), and
Propositions \ref{remark:whole-cone}, \ref{remark:index-1},
we obtain
$E_\sigma (\beta)=E_\sigma(\beta')$
for all facets.
By Lemma \ref{lemma:normal} (3),
the if-part follows from Theorem \ref{theorem:main}.
\qed

\begin{example}
Let 
$$
A=\left(
\begin{array}{cccc}
1 & 0 & 0 & 1\\
0 & 1 & 0 & 1\\
0 & 0 & 1 & -1
\end{array}
\right).
$$
Let $\beta\in {\bf Z}A={\bf Z}^d$.
Then by Theorem \ref{theorem:normal}, 
the $A$-hypergeometric system $M_A(\beta)$
is isomorphic to
$$
\begin{array}{rl}
M_A({}^t(0,0,0)) & \mbox{if}\quad 
\beta_1\geq 0, \beta_2\geq 0, \beta_1+\beta_3\geq 0, \beta_2+\beta_3\geq 0,\\
M_A({}^t(-1,0,1)) & \mbox{if}\quad
\beta_1< 0, \beta_2\geq 0, \beta_1+\beta_3\geq 0, \beta_2+\beta_3\geq 0,\\
M_A({}^t(0,-1,1)) & \mbox{if}\quad
\beta_1\geq 0, \beta_2< 0, \beta_1+\beta_3\geq 0, \beta_2+\beta_3\geq 0,\\
M_A({}^t(0,1,-1)) & \mbox{if}\quad
\beta_1\geq 0, \beta_2\geq 0, \beta_1+\beta_3< 0, \beta_2+\beta_3\geq 0,\\
M_A({}^t(1,0,-1)) & \mbox{if}\quad
\beta_1\geq 0, \beta_2\geq 0, \beta_1+\beta_3\geq 0, \beta_2+\beta_3< 0,\\
M_A({}^t(-1,-1,1)) & \mbox{if}\quad
\beta_1< 0, \beta_2< 0, \beta_1+\beta_3\geq 0, \beta_2+\beta_3\geq 0,\\
M_A({}^t(-1,0,0)) & \mbox{if}\quad
\beta_1< 0, \beta_2\geq 0, \beta_1+\beta_3< 0, \beta_2+\beta_3\geq 0,\\
M_A({}^t(0,-1,0)) & \mbox{if}\quad
\beta_1\geq 0, \beta_2< 0, \beta_1+\beta_3\geq 0, \beta_2+\beta_3< 0,\\
M_A({}^t(0,0,-1)) & \mbox{if}\quad
\beta_1\geq 0, \beta_2\geq 0, \beta_1+\beta_3< 0, \beta_2+\beta_3< 0,\\
M_A({}^t(-2,-1,1)) & \mbox{if}\quad
\beta_1< 0, \beta_2< 0, \beta_1+\beta_3< 0, \beta_2+\beta_3\geq 0,\\
M_A({}^t(-1,-2,1)) & \mbox{if}\quad
\beta_1< 0, \beta_2< 0, \beta_1+\beta_3\geq 0, \beta_2+\beta_3< 0,\\
M_A({}^t(-1,0,-1)) & \mbox{if}\quad
\beta_1< 0, \beta_2\geq 0, \beta_1+\beta_3< 0, \beta_2+\beta_3< 0,\\
M_A({}^t(0,-1,-1)) & \mbox{if}\quad
\beta_1\geq 0, \beta_2< 0, \beta_1+\beta_3< 0, \beta_2+\beta_3< 0,\\
M_A({}^t(-1,-1,0)) & \mbox{if}\quad
\beta_1< 0, \beta_2< 0, \beta_1+\beta_3< 0, \beta_2+\beta_3< 0.
\end{array}
$$
\end{example}

\section{Monomial curve case}

In this section, we conider $d=2$ case.
Let
$$ A \quad = \quad
\pmatrix{ 
1 & 1 & 1 & \cdots & 1 & 1 \cr
0 & i_2 & i_3 & \cdots & i_{n-1} & i_n \cr} $$ 
with $0 < \! i_2 \! < \! i_3 \! < \cdots < i_n $ relative prime integers.
Put $F_{\sigma_1}(s)=s_2$ and $F_{\sigma_2}(s)= i_n s_1 -s_2$.

We denote by ${\cal E}(A)$ the set of holes, i.e.,
\begin{eqnarray}
{\cal E}(A):&=&(({\bf N}A +{\bf Z}a_1)\cap ({\bf N}A +{\bf Z}a_n))
\setminus {\bf N}A
\\
&=&\{\, \beta\, |\,
E_{{\bf Q}_{\geq 0}A}(\beta)=\{ 0\},\,
E_{\sigma_1}(\beta)=\{ 0\},
\nonumber\\
&&\qquad\qquad\qquad
E_{\sigma_2}(\beta)=\{ 0\},\,
E_{\{ 0\}}(\beta)=\emptyset\,\}.
\end{eqnarray}

The rank of $M_A(\beta)$
is $d$ or $d+1$, and it equals $d+1$
if and only if $\beta\in {\cal E}(A)$
(see \cite{CDD}, \cite{sst-book}).

\begin{lemma}
\label{lemma:curve-index-1}
For any face $\tau$,
\begin{equation}
{\bf Z}A\cap ({\bf k}(A\cap\tau))
=
{\bf Z}(A\cap\tau).
\end{equation}
\end{lemma}

{\bf Proof.}
When $\tau$ is the whole cone ${\bf Q}_{\geq 0}A$ or the origin $\{ 0\}$,
the statement is trivial.

Note that $\beta$ belongs to ${\bf Z}A$ if and only if
$F_{\sigma_1}(\beta)\in {\bf Z}$,
$F_{\sigma_2}(\beta)\in {\bf Z}$,
and
$F_{\sigma_1}(\beta)+F_{\sigma_2}(\beta)
\in i_n{\bf Z}$.
Suppose $\beta\in
{\bf Z}A\cap ({\bf k}(A\cap\sigma_1))$.
Then $F_{\sigma_2}(\beta)\in i_n {\bf Z}$.
When $F_{\sigma_2}(\beta)= di_n $,
we have
$\beta =d a_n$.
\qed

\begin{corollary}
\begin{equation}
{\cal E}(A)= \{\,\beta \in {\bf Z}A\, |\,
F_{\sigma_1}(\beta)\in F_{\sigma_1}({\bf N}A),\,
F_{\sigma_2}(\beta)\in F_{\sigma_2}({\bf N}A)\,\}
\setminus {\bf N}A.
\end{equation}
\end{corollary}

{\bf Proof.}
This is immediate from Lemma \ref{lemma:curve-index-1}.
\qed

Theorem \ref{theorem:main} in the monomial curve case is as follows.

\begin{theorem}
Let $\beta,\beta'\in {\bf k}^d$.
\begin{enumerate}
\item
Suppose $\beta\notin {\cal E}(A)$.
Then $M_A(\beta')$ is isomorphic to $M_A(\beta)$
if and only if
$\beta -\beta'\in {\bf Z}A$, $\beta'\notin {\cal E}(A)$, and
$\{ \sigma_i : \, F_{\sigma_i} (\beta)\in F_{\sigma_i}({\bf N}A)\}=
\{ \sigma_i : \, F_{\sigma_i} (\beta')\in F_{\sigma_i}({\bf N}A)\}$.

\item
Suppose $\beta\in {\cal E}(A)$.
Then $M_A(\beta')$ is isomorphic to $M_A(\beta)$
if and only if
$\beta\in {\cal E}(A)$.
\end{enumerate}
\end{theorem}

{\bf Proof.}
(2) directly follows from Theorem \ref{theorem:main}.

The only-if-part of (1) follows from Theorem \ref{theorem:main}
by Proposition \ref{remark:whole-cone} (3).
Next suppose that $\beta -\beta'\in {\bf Z}A$, 
$\beta, \beta'\notin {\cal E}(A)$, and that
$\{ \sigma_i : \, F_{\sigma_i} (\beta)\in F_{\sigma_i}({\bf N}A)\}=
\{ \sigma_i : \, F_{\sigma_i} (\beta')\in F_{\sigma_i}({\bf N}A)\}$.
Then by Lemma \ref{lemma:curve-index-1},
Proposition \ref{remark:whole-cone} (3),
and Proposition \ref{remark:index-1} (2),
we have $E_{\sigma_i}(\beta)=E_{\sigma_i}(\beta')$ for $i=1,2$.
Moreover we know $E_{\{ 0\}}(\beta), E_{\{ 0\}}(\beta')=
\emptyset$ from Proposition \ref{remark:whole-cone} (2).
Hence $M_A(\beta)$ and $M_A(\beta')$ are isomorphic by Theorem
\ref{theorem:main}.
\qed

\begin{example}(\cite[Chapter 4]{sst-book})
Let
$$ A \quad = \quad
\pmatrix{ 
1 & 1 & 1  & 1 & 1\cr
0 & 2 & 4 & 7 & 9 \cr}. $$ 
Then
\begin{equation}
F_{\sigma_1}({\bf N}A)=\{\,
0, 2, 4, 6, 7, 8, 9,\ldots \,\},
\end{equation}
and
\begin{equation}
F_{\sigma_2}({\bf N}A)=\{\,
0, 2, 4, 5, 6, 7, 8, 9,\ldots \,\},
\end{equation}

Parameters in ${\bf Z}A={\bf Z}^2$
are decomposed into five parts
according to the isomorphism classes
of their corresponding $A$-hypergeometric systems:
\begin{enumerate}
\item
${\bf N}A$,
\item
$\{\, {}^t(\beta_1,\beta_2)\, |\,
\beta_2\in F_{\sigma_1} ({\bf N}A),\,
9\beta_1-\beta_2\notin F_{\sigma_2} ({\bf N}A)\,\}$,
\item
$\{\, {}^t(\beta_1,\beta_2)\, |\,
\beta_2\notin F_{\sigma_1} ({\bf N}A),\,
9\beta_1-\beta_2\in F_{\sigma_2} ({\bf N}A)\,\}$,
\item
$\{\, {}^t(\beta_1,\beta_2)\, |\,
\beta_2\notin F_{\sigma_1} ({\bf N}A),\,
9\beta_1-\beta_2\notin F_{\sigma_2} ({\bf N}A)\,\}$,
\item
${\cal E}(A)=\{ {}^t(2,10), {}^t(2,12), {}^t(3,19)\}$
: the set of holes.
\end{enumerate}
\end{example}

\section{Final remark}

Thanks to Theorem \ref{theorem:main},
all $D$-invariants of $A$-hypergeometric systems
can be described in terms of $E_\tau(\beta)$;
the characteristic cycles (in particular, the rank),
the monodromy representations, etc.
One of most recent results is given by
Tsushima (\cite{Tsushima})
on Laurent polynomial solutions.
He has proved that the vector space of Laurent polynomial solutions of
$M_A(\beta)$ has a basis consisting of canonical series
whose negative support corresponds to a face $\tau$ of ${\bf Q}_{\geq 0}A$
such that
$\dim \tau = |\{ \, a_j\, |\, a_j\in \tau\,\}|$,
and that $0\in E_\tau(\beta)$ but $0\notin E_{\tau'}(\beta)$ 
for any $\tau'\subset\tau$.
In particular, the dimension of the vector space of Laurent 
polynomial solutions
equals the cardinality of the set of such faces.
This is a generalization of the corresponding result
by Cattani, D'Andrea and Dickenstein (\cite{CDD})
in the monomial curve case.

\bigskip

Department of Mathematics

Hokkaido University

Sapporo, 060-0810

Japan

e-mail: saito@math.sci.hokudai.ac.jp


\begin{thebibliography}{ }

\bibitem{Adolphson}
Adolphson, A.(1994):
Hypergeometric functions and rings generated by monomials.
Duke Mathematical Journal {\bf 73}, 269--290.

\bibitem{CDD}
Cattani, E., D'Andrea, C., Dickenstein, A. (1998):
The ${\cal A}$-hypergeometric system associated
with a monomial curve.
Duke Mathematical Journal {\bf 99}, 179--207.

\bibitem{CDS}
Cattani, E., Dickenstein, A., Sturmfels, B. (1999):
Rational hypergeometric functions.
math.AG/9911030.

\bibitem{GKZ}
Gel'fand, I.M., Kapranov, M.M., Zelevinskii, A.V. (1990):
Generalized Euler integrals and $A$-hypergeometric functions.
Advances in Mathematics {\bf 84}, 255--271.

\bibitem{GZK-def}
Gel'fand, I.M.,  Zelevinskii, A.V., Kapranov, M.M. (1988):
Equations of hypergeometric type and Newton polyhedra.
Soviet Mathematics Doklady {\bf 37}, 678--683.


\bibitem{GZK}
Gel'fand, I.M.,  Zelevinskii, A.V., Kapranov, M.M. (1989):
Hypergeometric functions and toral manifolds.
Functional Analysis and its Applications {\bf 23}, 94--106.

\bibitem{Hosten-Thomas}
Hosten, S., Thomas, R.R. (1998):
Standard pairs and group relaxations in
integer programming,
to appear in Journal of Pure and Applied Algebra.

\bibitem{Oda}
Oda, T. (1988):
{\it Convex Bodies and Algebraic Geometry:
An Intriduction to the Theory of Toric Varieties.}
Ergebnisse der Mathematik und ihrer Grenzgebiete, 3.
Folge, Band 15, Springer, Heidelberg.


\bibitem{saito-normal}
Saito, M. (1994):
Normality of affine toric varieties associated with
Hermitian symmetric spaces.
Journal of Mathematical Society of Japan {\bf 46}, 699--724.

\bibitem{saito-fe}
Saito, M. (1995):
Contiguity Relations for the Lauricella Functions.
Funkcialaj Ekvacioj {\bf 38}, 37--58.

\bibitem{SymAlg}
Saito, M. (1996):
Symmetry algebras of normal ${\cal A}$-hypergeometric systems.
Hokkaido Mathematical Journal {\bf 25}, 591--619.

\bibitem{sst}
Saito, M., Sturmfels, B.,  Takayama, N. (1999):
Hypergeometric polynomials and integer programming.
Compositio Mathematica, {\bf 155}, 185--204.

\bibitem{sst-book}
Saito, M., Sturmfels, B.,  Takayama, N. (1999):
{\it Gr\"obner deformations of hypergeometric differential equations}.
Algorithms and Computation in Mathematics {\bf 6},
Springer, Berlin, Heidelberg, New York.

\bibitem{Saito-Takayama}
Saito, M., Takayama, N. (1994):
Restrictions of ${\cal A}$-hypergeometric systems and
connection formulas of the $\triangle_1\times\triangle_{n-1}$-hypergeometric 
function.
International Journal of Mathematics {\bf 5}, 537--560.

\bibitem{Tsushima}
Tsushima, T. :
Laurent polynomial solutions of $A$-hypergeometric systems.
in preparation.

\end{thebibliography}
\end{document}